\documentclass[final,11pt]{elsarticle}
\usepackage{amsmath}
\usepackage{amsthm}
\usepackage{amssymb}
\usepackage{mathtools}
\usepackage{booktabs}
\usepackage{bm}
\usepackage{xcolor}
\usepackage{tikz}

\usepackage{enumerate}
\usepackage{enumitem}
\setlist[enumerate,1]{label=(\arabic*),font=\textup,
leftmargin=7mm,labelsep=1.5mm,topsep=0mm,itemsep=-0.8mm}
\setlist[enumerate,2]{label=(\alph*),font=\textup,
leftmargin=7mm,labelsep=1.5mm,topsep=-0.8mm,itemsep=-0.8mm}

\usepackage{geometry}
\usepackage{microtype}

\usepackage[pdfstartview=FitH,bookmarksopen=false,
colorlinks=true,linkcolor=blue,citecolor=blue]{hyperref}

\newtheorem{theorem}{Theorem}[section]

\newtheorem{lemma}[theorem]{Lemma}
\theoremstyle{corollary}

\theoremstyle{observation}

\numberwithin{equation}{section}

\begin{document}
\begin{frontmatter}
\title{Nordhaus-Gaddum type inequality for the fractional matching number of a graph}

	\tnotetext[titlenote]{This work was supported by the National Nature Science Foundation
	of China (Nos. 11871040).}
	\author[]{Ting Yang}
    \ead{yangting_2019@shu.edu.cn}
    
  \author[]{Xiying Yuan\corref{correspondingauthor}}
   \cortext[correspondingauthor]{Corresponding author. Email: {\tt xiyingyuan@shu.edu.cn} (Xiying Yuan).}

\address{Department of Mathematics, Shanghai University, Shanghai 200444, P.R. China}


\begin{abstract}
The fractional matching number of a graph $G$, written as $\alpha^{\prime}(G)$, is the maximum size of a fractional matching of $G$. The following sharp lower bounds for a graph $G$ of order $n$ are proved, and all extremal graphs are characterized in this paper.\\
(1) $\alpha^{\prime}(G)+\alpha^{\prime}(\overline{G})\geq \frac{n}{2}$ for $n\geq2$.\\
(2) If $G$ and $\overline{G}$ are non-empty, then for $n\geq 28$, 
	$\alpha^{\prime}(G)+\alpha^{\prime}(\overline{G})\geq \frac{n+1}{2}$.\\
(3) If $G$ and $\overline{G}$ have no isolated vertices, then for $n\geq 28$,
$\alpha^{\prime}(G)+\alpha^{\prime}(\overline{G}) \geq \frac{n+4}{2}$.\\
\end{abstract}

\begin{keyword}
Nordhaus-Gaddum type inequality \sep
Fractional matching number \sep
Fractional Berge's theorem 
\end{keyword}
\end{frontmatter}

\section{Introduction}
Throughout this paper, all graphs are simple, undirected and finite. Undefined terminologies and notations can be found in \cite{3}. Let $G=(V(G), E(G))$ be a graph and $\overline{G}$ be its complement. $n$ will always denote the number of vertices of a given graph $G$. For a vertex $v\in V(G)$, its \textit{degree} $d_{G}(v)$ is the number of edges incident to it in $G$, its \textit{neighborhood}, denoted by $N(v)$, is the set of vertices, which are adjacent to $v$. 
An edge set $M$ of $G$ is called a \textit{matching} if any two edges in $M$ have no common vertices. 
The \textit{matching number} of a graph $G$, written $\alpha(G)$, is the number of edges in a maximum matching. As in \cite{2}, a \textit{fractional matching} of a graph $G$ is a function $f: E(G) \longrightarrow [0,1]$ such that $f(v)\leq 1$ for each vertex $v\in V(G)$, where $f(v)$ is the sum of $f(e)$ of edges incident to $v$. The \textit{fractional matching number} of $G$, written $\alpha^{\prime}(G)$, is the maximum value of $f(G)$ over all fractional matchings, where $f(G)$ denotes the sum of $f(e)$ of all edges in $G$. A \textit{fractional perfect matching} of a graph $G$ is a fractional matching $f$ with $f(G)=\alpha^{\prime}(G)=\frac{n}{2}$. Obviously, fractional matching is a generalization of matching. Choi et al. \cite{1} proved the difference and ratio of the fractional matching number and the matching number of graphs and characterized all infinite extremal family of graphs.\par
Nordhaus and Gaddum \cite{4} considered lower and upper bounds on the sum and the product of chromatic number $\chi(G)$ of a graph $G$ and its complement $\overline{G}$. They showed that 
\begin{equation*}
	2\sqrt{n}\leq \chi(G)+\chi(\overline{G})\leq n+1.
\end{equation*}
Since then, any bound on the sum or the product of an invariant in a graph $G$ and the same invariant in its complement $\overline{G}$ is called a \textit{Nordhaus-Gaddum type inequality}. 
The Nordhaus-Gaddum type inequality of various graph parameters has attracted much attention (see
\cite{6}\cite{7}\cite{8}\cite{9}\cite{10}\cite{11}\cite{12}). Aouchiche and Hansen \cite{13} wrote a stimulating survey on this topic, and we refer the reader to that article for additional information.
Chartrand and Schuster \cite{14} proved Nordhaus-Gaddum type result for the matching number of a graph. They showed that $\alpha(G)+\alpha(\overline{G})\geq \lfloor\frac{n}{2}\rfloor$. Laskar and Auerbach \cite{8} improved the bound by considering that $G$ and $\overline{G}$ contain no isolated vertices. They showed that $\alpha(G)+\alpha(\overline{G})\geq \lfloor\frac{n}{2}\rfloor+2$. Later, Lin et al. \cite{5} characterized all extremal graphs which attain the lower bounds of the results of Chartrand et al. and Laskar et al..
Motivated by them, we consider Nordhaus-Gaddum type inequality for the fractional matching number of a graph. First we prove some auxiliary results of fractional matching, which are selfcontained. Then we establish lower bounds on the sum of fractional matching number of a graph $G$ and its complement (see Theorem \ref{Theorem 3.3} and Theorem \ref{Theorem 3.5}). Moreover, we show those bounds are sharp. 
\section{Auxiliary results of fractional matching}

Based on an optimal fractional matching of graph $G$, we present a good partition of $V(G)$. We further characterize some properties of this partition.\par
\begin{lemma}\label{Lemma 2.1}(\cite{2})
	For any graph $G$, $2\alpha^{\prime}(G)$ is an integer. Moreover, there is a fractional matching $f$ for which 
	\begin{align*}
       f(G)=\alpha^{\prime}(G) 
	\end{align*}
	\noindent such that $f(e)\in \{0, \frac{1}{2}, 1\}$ for every edge $e$ of $G$. 
\end{lemma}

 An $f$ is called an \textit{optimal fractional matching} of graph $G$ in this paper, if we have \\
(1) $f(G)=\alpha^{\prime}(G)$.\\
(2) $f(e) \in \{0, \frac{1}{2}, 1\}$ for every edge $e$.\\
(3) $f$ has the greatest number of edges $e$ with $f(e)=1$. \\

In this paper, for a graph $G$, given a fractional matching $f$, an \textit{unweighted vertex} $v$ is a vertex with $f(v)=0$. A \textit{full vertex} $v$ is a vertex with $f(vw)=1$ for some edge $vw$ and we may call vertex $w$ is the \textit{full neighbour} of $v$. An \textit{$i-$edge} $e$ is an edge with $f(e)=i$. \textit{$\frac{1}{2}-$cycle} in a graph $G$ is an odd cycle induced by $\frac{1}{2}-$edges in $G$.  \\

\begin{lemma}\label{Lemma 2.2}
 Let $f$ be an optimal fractional matching of graph $G$. Then we have the following: \\
	$(1)$ (\cite{1}) The maximal subgraph induced by the $\frac{1}{2}-$edges is the union of odd cycles. \\
	$(2)$ (\cite{1}) The set of the unweighted vertices is an independent set of $G$. Furthermore, every unweighted vertex is adjacent only to a full vertex. \\
	$(3)$ No $\frac{1}{2}-$cycle has an unweighted vertex as a neighbour.
\end{lemma}
\noindent\textbf{Proof} By (2), every unweighted vertex is adjacent only to a full vertex, while every vertex on $\frac{1}{2}-$cycle is not a full vertex, then no $\frac{1}{2}-$cycle has an unweighted vertex as a neighbour. $\hfill\square $ \par


\begin{lemma}\label{Lemma 2.3}(\cite{2})
	Suppose that $f$ is a fractional matching of graph $G$. Then $f$ is a fractional perfect matching if and only if $f(v) =1$ for every vertex $v\in V(G)$. 
\end{lemma}
\begin{lemma}\label{Lemma 2.4}
		For any graph $G$, we have a partition $V(G)=V_{1}\dot{\cup}V_{2}$ with $|V_{1}|=2\alpha^{\prime}(G)$, and $G[V_{1}]$ contains a fractional perfect matching and $V_{2}$ is empty or an independent set. 
	\end{lemma}
	\noindent\textbf{Proof} Suppose $f$ is an optimal fractional matching of $G$. If $\alpha^{\prime}(G)=\frac{n}{2}$, we take $V_{1}=V(G)$. If $\alpha^{\prime}(G)<\frac{n}{2}$, set $V_{1}=\{v\in V(G)|f(v)>0\}$ and $V_{2}=V(G)\setminus V_{1}$. We will show that $f(v)=1$ for each vertex $v$ in $V_{1}$. Suppose to the contrary that there exists a vertex $v_{0}\in V_{1}$ with $0<f(v_{0})<1$, i.e., $f(v_{0})=\frac{1}{2}$, say $f(v_{0}v_{1})=\frac{1}{2}$ for some vertex $v_{1}\in V_{1}$. By Lemma \ref{Lemma 2.2} $(1)$, $v_{0}v_{1}$ lies in a $\frac{1}{2}-$cycle, then there exists a vertex $v_{t}( \neq v_{1})$ such that $f(v_{0}v_{t})=\frac{1}{2}$. Thus $f(v_{0})=1$, which is a contradiction. Thus $G[V_{1}]$ contains a fractional perfect matching by Lemma \ref{Lemma 2.3}. Since $V_{2}$ is a set of the unweighted vertices, $V_{2}$ is an independent set by Lemma \ref{Lemma 2.2} $(2)$. Furthermore, we have $\alpha^{\prime}(G)=\alpha^{\prime}(G[V_{1}])=\frac{|V_{1}|}{2}$.  $\hfill\square $\par
	\vspace{0.5cm}
	For a graph $G$ with $\alpha^{\prime}(G)=t$, we write $V(G)=V_{1}\dot\cup V_{2}$ according to the results of Lemma \ref{Lemma 2.4}, and let $s$ be the maximum number of independent edges in $[V_{1}, V_{2}]$. Now we will further decompose $V_{1}=V_{11}\dot\cup V_{12}$ and
	$V_{2} = V_{21} \dot\cup V_{22}$ such that $[V_{11}, V_{21}]$ contains exactly $s$ independent edges. Obviously, $|V_{11}| =|V_{21}| = s$ holds. We call $V(G)=(V_{11}\cup V_{12})\cup(V_{21}\cup V_{22})$ a \textit{good partition} of $G$ in this paper (see Figure 1).
\begin{lemma}\label{Lemma 2.5}
	Let $V(G)=(V_{11}\cup V_{12})\cup(V_{21}\cup V_{22})$ be a good partition of $G$ with $|V_{11}|=s$ and $f$ be the corresponding optimal fractional matching of $G$.\\
	$(1)$ If $f(uv)=1$ for some edge $uv$ in $G$, then there is no vertex in $V_{2}\cap N(u)\cap 
	N(v)$. \\
	$(2)$ If $e$ is an edge in $G[V_{11}]$, then $f(e)=0$. \\
	$(3)$ Each vertex in $V_{11}$ is a full vertex.  
\end{lemma}
\begin{figure}
	\tikzstyle{format}=[draw,thin,fill=white]
	\begin{center}
		\begin{tikzpicture}
		\draw (1,1) ellipse (3cm and 1.5cm)node[right=3cm]{$V_2$}node[below=2.0cm]{$G$}node[below=3.0cm]{Figure 1: A good partition of graph $G$};
		\draw (1,5) ellipse (3cm and 1.5cm)node[right=3cm]{$V_1$};
		\draw (-0.5,1) ellipse (1.0cm and 0.6cm)node[below=0.6cm]{$V_{21}$};
		\draw (2.5,1) ellipse (1.4cm and 0.9cm)node[below=0.78cm]{$V_{22}$};
		\draw (-0.5,5) ellipse (1.0cm and 0.6cm)node[above=0.6cm]{$V_{11}$};
		\draw (2.5,5) ellipse (1.4cm and 0.9cm)node[above=0.78cm]{$V_{12}$};
		\draw (2.5,5) ellipse (1.0cm and 0.6cm)node[below right=0.10cm]{$X$};
		
		\filldraw (-1.2,5.3) circle (1.2pt);
		\filldraw (0.2,4.7) circle (1.2pt);
		\filldraw (-0.85,5.15) circle (0.5pt)node[below right= 0.2cm]{$s$};
		\filldraw (-0.5,5) circle (0.5pt);
		\filldraw (-0.15,4.85) circle (0.5pt);
		
		\filldraw (-1.2,1) circle (1.2pt);
		\filldraw (0.2,1) circle (1.2pt);
		\filldraw (-0.85,1) circle (0.5pt);
		\filldraw (-0.5,1) circle (0.5pt);
		\filldraw (-0.15,1) circle (0.5pt);
		
		\filldraw (2.5,5.3) circle (1.2pt);
		\filldraw (2.5,4.7) circle (1.2pt);
		\filldraw (2.5,5.15) circle (0.5pt);
		\filldraw (2.5,5) circle (0.5pt);
		\filldraw (2.5,4.85) circle (0.5pt);
		
		\draw (-1.2,5.3) -- (-1.2,1);
		\draw (0.2,4.7) -- (0.2,1);
		\draw (-1.2,5.3) -- (2.5,5.3);
		\draw (0.2,4.7) -- (2.5,4.7);
		
		\filldraw (1.4,1) circle (1.2pt);
		\filldraw (3.5,1) circle (1.2pt);
		\filldraw (1.925,1) circle (0.5pt);
		\filldraw (2.45,1) circle (0.5pt);
		\filldraw (2.975,1) circle (0.5pt);
		\filldraw (3.65,5) circle (0.5pt);
		\filldraw (3.65,5.3) circle (1.2pt);
		\filldraw (3.65,4.7) circle (1.2pt);
		
		\end{tikzpicture}
	\end{center}
\end{figure}\par

\noindent\textbf{Proof} (1) Suppose there exists a vertex $w\in V_{2}\cap N(u)\cap 
N(v)$. Since $w$ is an unweighted vertex, we have $f(uw)=f(vw)=0$. Now set $f^{*}(uw)=f^{*}(vw)=f^{*}(uv)=\frac{1}{2}$ and other assignments remain unchanged. Then $f^{*}(G)=\alpha^{\prime}(G)+\frac{1}{2}$, which is a contradiction. \par
(2) Suppose $uv$ is an edge in $G[V_{11}]$ with $f(uv)=a>0$. Let $ux$ and $vy$ be two independent edges in $[V_{11}, V_{21}]$. Now set $f^{*}(ux)=f^{*}(vy)=1$ and $f^{*}(uz)=f^{*}(vw)=0$ for any $z\in N(u)\setminus\{x\}$ and $w\in N(v)\setminus\{y\}$ and other assignments remain unchanged. Then $f^{*}(G)-f(G)\geq a$, which is a contradiction to the choice of $f$. \par
(3) Suppose $u\in V_{11}$ is not a full vertex. Then there exists a vertex $v$ such that $f(uv)=\frac{1}{2}$. Furthermore, $uv$ lies in a $\frac{1}{2}-$cycle by Lemma \ref{Lemma 2.2} (1). While $u$ has a neighbour in $V_{21}$, which contradicts to Lemma \ref{Lemma 2.2} (3). $\hfill\square $ \\

From Lemma \ref{Lemma 2.5} (3), we know that each vertex in $V_{11}$ is a full vertex. Denote by $X$ the set of all full neighbours of vertices in $V_{11}$, and then $X \subseteq V_{12}$ by Lemma \ref{Lemma 2.5} (2). It is obvious that $|X|=|V_{11}|=s$. Furthermore, we have the following results. 
\begin{lemma}\label{Lemma 2.6}
	$(1)$ If $|X|=s\geq 2$, then $X$ is an independent set of $G$.\\
	$(2)$ There is no edge between $V_{2}$ and $X$ in $G$.
\end{lemma}
\noindent\textbf{Proof} (1) Suppose there exist two vertices $u, v \in X$ such that $uv\in E(G)$ and $f(gu)=f(hv)=1$ for  vertices $g, h$ in $V_{11}$. Let $gw$ and $hy$ be two independent edges in $[V_{11}, V_{21}]$. Now set $f^{*}(gu)=f^{*}(hv)=0$ and $f^{*}(uv)=f^{*}(gw)=f^{*}(hy)=1$ and other assignments remain unchanged. Then $f^{*}(G)=\alpha^{\prime}(G)+1$, which is a contradiction.  \par
(2) Suppose to the contrary that there exists an edge $vx$ with $v\in X$ and $x\in V_{21}$. Then there exists a vertex $u\in V_{11}$ such that $uv \in E(G)$ and $f(uv)=1$. We have $ux\notin E(G)$ in virtue of Lemma \ref{Lemma 2.5} (1). Then there exists a vertex say $y\in V_{21}$ such that $uy$ is one of the independent edges in $[V_{11}, V_{21}]$. Now we may modify $f$ to $f^{*}$. Let $f^{*}(uv)=0, f^{*}(vx)=f^{*}(uy)=1$ and other assignments remain unchanged. Then $f^{*}(G)=\alpha^{\prime}(G)+1$, which is a contradiction. It is obvious that there is no edge between $V_{22}$ and $X$ in $G$. Thus there is no edge between $V_{2}$ and $X$ in $G$.   $\hfill\square $ \par
\vspace{0.5cm}

\section{Graphs with small fractional matching number}

In virtue of the fractional Berge$^{\prime}$s theorem, we will characterize some graphs with small fractional matching number in this section.
\begin{lemma}\label{Lemma 2.7}(\cite{2})
	For any graph $G$ of order $n$, we have
	\begin{equation*}
	\alpha^{\prime}(G)=\frac{1}{2} \bigg(n-\mathop{max}\limits_{S\subseteq V(G)}\Big\{i\big (G-S\big)-|S|\Big\}\bigg)		,
	\end{equation*}
	where $i(G-S)$ denotes the number of isolated vertices of $G-S$.
\end{lemma}
By the fractional Berge$^{\prime}$s theorem, we immediately have the following results. 
\begin{lemma}\label{Lemma 2.8}
	$(1)$ For a graph $G$ of order $n$, $\alpha ^{\prime}(G)=1$ if and only if $G \cong K_{1, k} \cup (n-1-k) K_{1}$, where $k\geq 1$. \\
	$(2)$ For a graph $G$ of order $n$, $\alpha ^{\prime} (G)=\frac{3}{2}$ if and only if $G \cong C_{3} \cup (n-3) K_{1}$.
\end{lemma}
	
		
		

Let $K_{2}(p, q; \ell)(p\geq q)$ be the graph obtained by attaching $p$ pendent edges at one vertex of $K_{2}$ called $uv$, $q$ pendent edges at the other vertex of $K_{2}$ and having $\ell$ vertices in $N(u)\cap N(v)$. 
\begin{lemma}\label{Lemma 2.9}
	For a graph $G$ of order $n$, $\alpha^{\prime}(G)=2$ if and only if one of the following situations occurs:  \\
	$(1)$ $2K_{2}\cup (n-4)K_{1}\subseteq G \subseteq K_{4} \cup (n-4)K_{1}$. \\
	$(2)$ $2K_{2}\cup (n-4)K_{1}\subseteq G\subseteq K_{2}(0, 0; n-2)$. 
\end{lemma}
\noindent\textbf{Proof} By considering the fact that if $G_{1}\subseteq G_{2}$, then $\alpha^{\prime}(G_{1}) \leq \alpha^{\prime}(G_{2})$ and $\alpha^{\prime}(K_{4}\cup(n-4)K_{1})=\alpha^{\prime}(K_{2}(0, 0; n-2))=\alpha^{\prime}(2K_{2}\cup(n-4)K_{1})=2$, the sufficiency part is correct. \par 
To show the necessity part, by Lemma \ref{Lemma 2.7}, suppose $S \subseteq V(G)$ such that $n-4 = i(G-S)-|S|$. Since $ i(G-S) \leq n-|S|$, it follows that $|S|\leq 2$.  \par
If $|S| = 0$, then $i(G)=n-4$. $\alpha^{\prime}(G)=2$ implies $G$ contains two independent edges. Thus we have $2K_{2}\cup (n-4)K_{1}\subseteq G \subseteq K_{4} \cup (n- 4)K_{1}$.  \par
If $|S| =1$, then $i(G-S)=n-3$, and $G-S$ contains a subgraph $F=K_{2}$. $\alpha^{\prime}(G)=2$ implies that there is an edge between $S$ and $(V(G)\setminus\{S\cup V(F)\})$. Thus $2K_{2}\cup(n-4)K_{1}$ is a subgraph of $G$. Since there are at most $n-3$ edges in $[S, V(G)\setminus\{S\cup V(F)\}]$ and there are at most two edges in $[S, V(F)]$, $G$ is a subgraph of $K_{2}(n-3, 0; 1)$. Thus $2K_{2}\cup(n-4)K_{1}\subseteq G \subseteq K_{2}(n-3, 0; 1)$.\par
If $|S| = 2$, then $i(G-S)=n-2$. $\alpha^{\prime}(G)=2$ implies $G$ contains two independent edges in $[S, V(G)\setminus S]$. Thus we have $2K_{2}\cup (n-4)K_{1}\subseteq G$. Since there are at most $2(n-2)$ edges in $[S, V(G)\setminus S]$ and there is at most one edge in $S$, $G$ is a subgraph of $K_{2}(0, 0; n-2)$. Thus $2K_{2}\cup (n-4)K_{1}\subseteq G \subseteq K_{2}(0, 0; n-2)$.   
 Noting that $K_{2}(n-3, 0; 1)\subseteq K_{2}(0, 0; n-2)$, we have $2K_{2}\cup (n-4)K_{1}\subseteq G\subseteq K_{2}(0, 0; n-2)$. $\hfill\square $\\
\begin{lemma}\label{Lemma 2.10}
	Let $H$ be the graph of order $n$ obtained by attaching $n-4$ pendent edges at one vertex of $K_{4}$. For a graph $G$ of order $n$, $\alpha^{\prime}(G)=\frac{5}{2}$ if and only if one of the following situations occurs:  \\
	$(1)$ $C_{5}\cup (n-5)K_{1}\subseteq G \subseteq K_{5} \cup (n-5)K_{1}$. \\
$(2)$ $C_{3}\cup K_{2}\cup (n-5)K_{1}\subseteq G \subseteq K_{5}\cup (n-5)K_{1}$. \\ $(3)$ $C_{3}\cup K_{2}\cup (n-5)K_{1} \subseteq G\subseteq H$.
\end{lemma}
\noindent\textbf{Proof} The sufficiency part is obvious. To show the necessity part, by Lemma \ref{Lemma 2.7}, suppose $S \subseteq V(G)$ such
that $n-5 = i(G-S)-|S|$. Since $ i(G-S) \leq n-|S|$, it follows that $|S|\leq\frac{5}{2}$. Furthermore, $|S|=2$ does not occur. Otherwise we have $i(G-S)=n-3$. While we have $G-S=(n-3)K_{1} \cup K_{1}$, which implies that $i(G-S)=(n-2)$.\par
If $|S| = 0$, then $i(G)=n-5$. $\alpha^{\prime}(G)=\frac{5}{2}$ implies $G$ contains an odd cycle. If $G$ contains $C_{3}$ as a subgraph, then we have $C_{3}\cup K_{2}\cup (n-5)K_{1}\subseteq G \subseteq K_{5}\cup (n-5)K_{1}$. If $G$ contains $C_{5}$ as a subgraph, then we have $C_{5}\cup (n-5)K_{1}\subseteq G \subseteq K_{5}\cup (n-5)K_{1}$. \par
If $|S| = 1$, then $i(G-S)=n-4$. $\alpha^{\prime}(G)=\frac{5}{2}$ implies that $G-S$ contains $C_{3}$ as a subgraph, and there exists at least one edge in $[S, V(G)\setminus\{S\cup V(C_{3})\}]$. Thus $C_{3}\cup K_{2}\cup (n-5)K_{1}$ is a subgraph of $G$. Furthermore, we have $G\subseteq H$.$\hfill\square $\\

Based on the results of Lemma \ref{Lemma 2.8}, Lemma \ref{Lemma 2.9} and Lemma \ref{Lemma 2.10}, we have the following results.
\begin{lemma}\label{Lemma 2.11}
	Let $G$ be a graph of order $n$. Then the following statements hold. \\
	$(1)$ If $\alpha^{\prime}(G) =1$ and $n\geq 4$, then $\alpha^{\prime}(G)+\alpha^{\prime}(\overline{G})\geq\frac{n+1}{2}$.\\
	$(2)$ If $\alpha^{\prime}(G) =\frac{3}{2}$ and $n\geq 6$, then $\alpha^{\prime}(G)+\alpha^{\prime}(\overline{G})=\frac{n+3}{2}$. \\
	$(3)$ If $\alpha^{\prime}(G) =2$ and $G$ is not isomorphic to $K_{2}(0, 0; \ell)$ for $\ell \geq 2$, then 
	$\alpha^{\prime}(G)+\alpha^{\prime}(\overline{G})\geq \frac{n+3}{2}$ for $n\geq 8$.\\
	$(4)$ If $\alpha^{\prime}(G) =2$ and $G$ is isomorphic to $K_{2}(0, 0; \ell)$ for $\ell \geq 2$, then $\alpha^{\prime}(G)+\alpha^{\prime}(\overline{G})=\frac{n+2}{2}$.\\
	$(5)$ If $\alpha^{\prime}(G) =\frac{5}{2}$ and $n\geq 7$, then $\alpha^{\prime}(G)+\alpha^{\prime}(\overline{G})\geq\frac{n}{2}+2$.\\ 
	Moreover, the equalities in $(1)$, $(3)$ and $(5)$ hold if and only if $G$ contains exactly one vertex with degree $n-1$.
\end{lemma}

\section{Nordhaus-Gaddum-type bounds for the fractional matching number }
In this section, we will prove the Nordhaus-Gaddum-type bounds for the fractional matching number (see Theorem \ref{Theorem 3.3} and Theorem \ref{Theorem 3.5}).\par
In a graph $G$ with $\alpha^{\prime}(G)=t$, we try to find a collection $E$ of some independent edges in $\overline{G}$ and assign $1$ to each edge of $E$, and then find a long cycle $C$ (having no common vertices with $E$) and assign $\frac{1}{2}$ to each edge of $C$. By this way, we get a lower bound of $\alpha^{\prime}(\overline{G})$.

\begin{lemma}\label{Lemma 3.1}
	Let $V(G)=(V_{11}\cup V_{12})\cup (V_{21}\cup V_{22})$ be a good partition of graph $G$ of order $n$ with $\alpha^{\prime}(G)=t \leq \frac{n}{4}$ and $|V_{11}|=s$. Then the following statements hold.\\
		$(1)$ $\alpha^{\prime}( \overline{G} )\geq \frac{n-s}{2}\geq \frac{n-t}{2}$. \\
		$(2)$ If both $G$ and $\overline{G}$ contain no isolated vertices and $s\geq 1$, then $\alpha^{\prime}(\overline{G}) \geq \frac{n-s+1}{2}$.\\  
		$(3)$ If both $G$ and $\overline{G}$ contain no isolated vertices and $s=t\geq 3$, then $\alpha^{\prime}(\overline{G}) \geq \frac{n-s+2}{2}$.  
\end{lemma}
\noindent \textbf{Proof} Since $V(G)=(V_{11} \cup V_{12})\cup (V_{21}\cup V_{22})$ is a good partition, $\overline{G}$[$V_{2}$] is a clique, and each vertex in $V_{12}$ is adjacent to each vertex in $V_{22}$ in $\overline{G}$. The assumption $n\geq 4t$ insures $|V_{22}|=n-2t-s\geq 2t-s$. 

(1) If $n-2t-(2t-s)=0$, that is $s=0$ and $|V_{12}|=|V_{22}|=\frac{n}{2}=2t$, then we may assign number $1$ to $2t$ independent edges in $\overline{G}[V_{12}, V_{22}]$ for a fractional matching of $\overline{G}$. We have $\alpha^{\prime}(\overline{G})\geq 2t=\frac{n}{2}$. \par
If $n-2t-(2t-s)=[(n-2t-s)-(2t-s)]+s=1$, then $s\leq 1$. There exist $u\in V_{12}$ and $v, w \in V_{2}$ forming a cycle for a fractional matching of $\overline{G}$. We may assign number $1$ to $2t-s-1$ independent edges in $\overline{G}[V_{12}\setminus\{u\}, V_{2}\setminus\{v, w\}]$, and $\frac{1}{2}$ to edges of $C_{3}$ induced by vertices $u, v$ and $w$. Thus we have $\alpha^{\prime}(\overline{G})\geq 2t-s-1+\frac{3}{2}=\frac{n-s}{2}$. \par
If $n-2t-(2t-s)=2$, then we may choose vertices $u, v \in V_{2}$ and assign number $1$ to edge $uv$ in $\overline{G}$ and $2t-s$ independent edges in $\overline{G}[V_{12}, V_{2}\setminus\{u, v\}]$ for a fractional matching of $\overline{G}$. Thus we have $\alpha^{\prime}(\overline{G})\geq 2t-s+1=\frac{n-s}{2}$.\par
When $n-2t-(2t-s)\geq 3$, we may assign number $1$ to $2t-s$ independent edges in $\overline{G}[V_{12}, V_{22}]$, and $\frac{1}{2}$ to each edge of a cycle with length $n-2t-(2t-s)$ in $\overline{G}[V_{2}]$ for a fractional matching of $\overline{G}$. Then we have
\begin{align*}
\alpha^{\prime}(\overline{G}) \geq 2t-s+\frac{n-2t-(2t-s)}{2} =\frac{n-s}{2}.
\end{align*}\par
Since $s\leq t$, we have $\alpha^{\prime}( \overline{G} )\geq \frac{n-s}{2}\geq \frac{n-t}{2}$.\par

$(2)$ Now we suppose $s\geq1$ and $\overline{G}$ contains no isolated vertices. Let $u$ be a vertex in $V_{11}$ and $v$ be one of its neighbours in $\overline{G}$. First we suppose $v\in V_{12}$. Now we define a fractional matching of $\overline{G}$. Let $\bar{f}(uv)=1$ in $\overline{G}$. Noting that $n-2t-(2t-s-1)\geq 2$. If $n-2t-(2t-s-1)=2$, that is $s=1$ and $|V_{12}|=|V_{22}|=2t-1=\frac{n}{2}-1$, then we may choose vertices $v _{1}\in V_{22}$ and $v_{2}\in V_{21}$ and assign number $1$ to edge $v_{1}v_{2}$ in $\overline{G}$ and $2t-s-1$ independent edges in $\overline{G}[V_{12}\setminus\{v\}, V_{22}\setminus\{v_{1}\}]$. Thus we have $\alpha^{\prime}(\overline{G})\geq 2+2t-s-1=\frac{n-s+1}{2}$.\par
When $n-2t-(2t-s-1)\geq 3$, we may assign number $1$  to  $2t-s-1$ independent edges in $\overline{G}[V_{12}\setminus \lbrace v \rbrace, V_{22}]$, and $\frac{1}{2}$ to each edge of a cycle with length $n-2t-(2t-s-1)$ in $\overline{G}[V_{2}]$ for a fractional matching of $\overline{G}$. Then we have
\begin{align*}
\alpha^{\prime}(\overline{G}) \geq 1+(2t-s-1)+\frac{n-2t-(2t-s-1)}{2} =\frac{n-s+1}{2}.
\end{align*}\par


Now suppose $v$ lies in $V_{11}$ or $V_{2}$. If $v\in V_{11}$ or $v\in V_{21}$, then $s\geq 2$. When $v\in V_{22}$, noting $G$ contains no isolated vertices, all possible neighbours of $v$ are in $V_{11}$ and $uv\notin E(G)$, then $|V_{11}|=s\geq 2$. Hence, we have $|X|=|V_{11}|\geq 2$ and there is an edge $v_{1}v_{2}$ in $\overline{G}[X]$ by Lemma \ref{Lemma 2.6} (1). Noting that  $n-2t-(2t-s-2)-1\geq 1+s\geq 3$. For a fractional matching of $\overline{G}$, we may assign number $1$ to  the edges $v_{1}v_{2}$, $uv$ and $2t-s-2$ independent edges in $\overline{G}[V_{12}\setminus \lbrace v_{1}, v_{2} \rbrace, V_{22}\setminus\{v\}]$, and $\frac{1}{2}$ to each edge of a cycle with length $n-2t-(2t-s-2)-1$ in $\overline{G}[V_{2}]$.
Thus we have 
\begin{align*}
\alpha^{\prime}(\overline{G}) \geq 1+1+(2t-s-2)+\frac{n-2t-(2t-s-2)-1}{2} =\frac{n-s+1}{2}.
\end{align*}\par
(3) Now we suppose that $s=t\geq 3$ and $\overline{G}$ contains no isolated vertices. The assumption $s=t$ implies $X=V_{12}$. By Lemma \ref{Lemma 2.6}, we obtain that  $\overline{G} [V_{12}]=K_{s}$ and each vertex in $V_{12}$ is adjacent to each vertex in $V_{2}$ in $\overline{G}$. Let $u$ be a vertex in $V_{11}$ and $v$ be one of its neighbours in $\overline{G}$. Noting that $n-2t-(2t-s)\geq s\geq 3$. Now we define a fractional matching $\bar{f}$ of $\overline{G}$. We suppose $v \in V_{11}$. First let $\bar{f}(uv)=1$ in $\overline{G}$ and we may assign number $1$ to $2t-s$ independent edges in $\overline{G}[V_{12}, V_{22}]$, and $\frac{1}{2}$ to each edge of a cycle with length $n-2t-(2t-s)$ in $\overline{G}[V_{2}]$. Thus we have 
\begin{align*}
\alpha^{\prime}(\overline{G}) \geq 1+2t-s+\frac{n-2t-(2t-s)}{2} =\frac{n-s+2}{2}.
\end{align*}\par
If $v\in V_{2}$, then there exists a vertex say $w \neq u \in V_{11}$ such that $wv\in E(G)$. Since $d_{G} (w) \leq n-2$, there exists a vertex $w^{\prime}$ such that $ww^{\prime}\notin E(G)$. If $w^{\prime}\in V_{11}$, it is the case discussed above. If $w^{\prime}\in V_{12}$, then let $\bar{f}(uv)=\bar{f}(ww^{\prime})=1$ in $\overline{G}$ and we may assign number $1$ to $2t-s-1$ independent edges in $\overline{G}[V_{12}\setminus\{w^{\prime}\}, V_{2}\setminus\{v\}]$, and $\frac{1}{2}$ to each edge of a cycle with length $n-2t-(2t-s-1)-1$ in $\overline{G}[V_{2}]$. Thus we have
\begin{align*}
\alpha^{\prime}(\overline{G}) \geq 2+(2t-s-1)+\frac{n-2t-(2t-s-1)-1}{2} =\frac{n-s+2}{2}.
\end{align*}\par
If $w^{\prime}\in V_{2}$, then let $\bar{f}(uv)=\bar{f}(ww^{\prime})=1$ in $\overline{G}$. Since $\overline{G}[V_{12}]=K_{s}$, there exist $x_{1}, x_{2} \in V_{12}$ such that $x_{1}x_{2} \in E(\overline{G})$. Let $\bar{f}(x_{1}x_{2})=1$ in $\overline{G}$ and we may assign number $1$ to $2t-s-2$ independent edges in $\overline{G}[V_{12}\setminus\{x_{1}, x_{2}\}, V_{2}\setminus\{v, w^{\prime}\}]$, and $\frac{1}{2}$ to each edge of a cycle with length $n-2t-2-(2t-s-2)$ in $\overline{G}[V_{2}]$. Thus we have 
\begin{align*}
\alpha^{\prime}(\overline{G}) \geq 3+2t-s-2+\frac{n-2t-2-(2t-s-2)}{2} =\frac{n-s+2}{2}.
\end{align*}\par
When $v\in V_{11}$ or $v\in V_{2}$, we obtain the desired inequality. So in the following, we may assume that $G[V_{11}]$ is a clique and each vertex in $V_{11}$ is adjacent to each vertex in $V_{2}$ in $G$. We suppose $v$ lies in $V_{12}$. From the definition of set $X (=V_{12})$, we have $v^{\prime}v\in E(G)$ for some $v^{\prime} \in V_{11}$. Since $v^{\prime}$ is not an isolated vertex of $\overline{G}$, we may assume $v^{\prime}v^{\prime\prime}\in E(\overline{G})$ for some vertex $v^{\prime\prime} \in V_{12}$. Now we construct a fractional matching $\bar{f}$ of $\overline{G}$ with $\bar{f}(uv)=\bar{f}(v^{\prime}v^{\prime\prime})=1$ in $\overline{G}$ and we may assign number $1$ to $2t-s-2$ independent edges in $[V_{12}\setminus \lbrace v, v^{\prime\prime} \rbrace, V_{22}]$ in $\overline{G}$, and  $\frac{1}{2}$ to each edge of a cycle with length $n-2t-(2t-s-2)$ in $\overline{G}[V_{2}]$. Thus we have 
\begin{align*}
\alpha^{\prime}(\overline{G}) \geq 2+2t-s-2+\frac{n-2t-(2t-s-2)}{2} =\frac{n-s+2}{2}.
\end{align*}\par
This completes the proof.       $\hfill\square $\\

\begin{lemma}\label{Lemma 3.2}
	$(1)$ If $\alpha^{\prime}(G)=\lfloor\frac{n}{4}\rfloor+\frac{1}{2}$ when $n\equiv 0, 1(mod4)$ or  $\alpha^{\prime}(G)=\lfloor\frac{n}{4}\rfloor+\frac{3}{2}$ when $n\equiv 2, 3(mod4)$, then $\alpha^{\prime}(\overline{G} )\geq \frac{n}{4}+3$ for $n\geq 28$. \\
	$(2)$ If $\alpha^{\prime}(G)=\lfloor\frac{n}{4}\rfloor+1$, then $\alpha^{\prime}(\overline{G} )\geq \frac{n}{4}+3$ for $n\geq 28$. 
\end{lemma}
\noindent \textbf{Proof} Recall $V(G)=(V_{11}\cup V_{12}) \cup (V_{21}\cup V_{22})$, $|V_{11}|=s$ and $\alpha^{\prime}(G)=t$. We will prove $\alpha^{\prime}(\overline{G})\geq \frac{n-t}{2}$. When $s\leq1$, noting that $|V_{1}|=2\alpha^{\prime}(G)>\frac{n}{2}$ and $|V_{12}|=2t-s> n-2t-s=|V_{22}|$. Then we may assign number $1$ to $n-2t-s$ independent edges in $\overline{G}[V_{12}, V_{22}]$ for a fractional matching of graph $\overline{G}$. Thus
\begin{align*}
\alpha^{\prime}(\overline{G}) \geq n-2t-s\geq \frac{n-t}{2}
\end{align*} for $n \geq 20$.\par
When $s\geq 2$, there exist two vertices $v_{1}, v_{2}\in X$ such that $v_{1}v_{2}\notin E(G)$ by Lemma \ref{Lemma 2.6} (1). We will define a fractional matching $\bar{f}$ of $\overline{G}$. Noting that $|V_{22}|\geq 3$. \par

(1) If $\alpha^{\prime}(G)=\lfloor\frac{n}{4}\rfloor+\frac{1}{2}$ when $n\equiv 0, 1(mod4)$ or $\alpha^{\prime}(G)=\lfloor\frac{n}{4}\rfloor+\frac{3}{2}$ when $n\equiv 2, 3(mod4)$, then $s<t$ and there exists a $\frac{1}{2}-$cycle in $V_{12}\setminus X$. By Lemma \ref{Lemma 2.2} (3), from the $\frac{1}{2}-$cycle of $G$, we may take two different vertices $v^{\prime}, v^{\prime\prime}$. Let $\bar{f}(v^{\prime}x^{\prime})=\bar{f}(v^{\prime\prime}x^{\prime\prime})=1$ in $\overline{G}$ for $\{x^{\prime}, x^{\prime\prime}\}\subseteq V_{21}$. There are $s-2$ independent edges in $\overline{G}[X\setminus\{v_{1}, v_{2}\}, V_{21}\setminus\{x^{\prime}, x^{\prime\prime}\}]$ by Lemma \ref{Lemma 2.6} (2) and we may assign number $1$ to them. Noting that $|V_{1}|-s-2 \leq |V_{2}|$, we may assign number $1$ to $2t-2s-2$ independent edges in $\overline{G}[V_{12}\setminus\{X, v^{\prime}, v^{\prime\prime}\}, V_{22}]$. When $n-2t-(2t-s-2)=0$, we may assign number $1$ to edge $v_{1}v_{2}$ in $\overline{G}$. Thus $\alpha^{\prime}(\overline{G}) \geq 1+2t-s-2=2t-s-1=\frac{n-s}{2}>\frac{n-t}{2}$.\par
If $n-2t-(2t-s-2)=1$, that is $n=4t-s-1$ and $|V_{22}|=2t-2s-1$, then we may choose $w\in V_{22}$ and assign number $\frac{1}{2}$ to edges of $C_{3}$ induced by $v_{1}, v_{2}$ and $w$ in $\overline{G}$, and $1$ to $2t-2s-2$ independent edges in $\overline{G}[V_{12}\setminus\{X, v^{\prime}, v^{\prime\prime}\}, V_{22}]$. Thus we have $\alpha^{\prime}(\overline{G}) \geq \frac{3}{2}+2t-s-2=2t-s-\frac{1}{2}=\frac{n-s}{2}>\frac{n-t}{2}$.\par
If $n-2t-(2t-s-2)=2$, then we may choose $w, w^{\prime}\in V_{22}$ and assign number $1$ to $2t-2s-2$ independent edges in $\overline{G}[V_{12}\setminus\{X, v^{\prime}, v^{\prime\prime}\}, V_{22}\setminus\{w, w^{\prime}\}]$. Let $\bar{f}(v_{1}v_{2})=\bar{f}(ww^{\prime})=1$ in $\overline{G}$. Thus we have $\alpha^{\prime}(\overline{G}) \geq 2+2t-s-2=2t-s=\frac{n-s}{2}>\frac{n-t}{2}$.\par
When $n-2t-(2t-s-2)\geq 3$, we may assign number $\frac{1}{2}$ to each edge of a cycle with length $n-2t-(2t-s-2)$ in $\overline{G}[V_{22}]$. Let $\bar{f}(v_{1}v_{2})=1$ in $\overline{G}$. Thus
\begin{align*}
\alpha^{\prime}(\overline{G}) \geq 1+2t-s-2+\frac{n-2t-(2t-s-2)}{2} =\frac{n-s}{2}> \frac{n-t}{2}.
\end{align*}\par

(2) When $\alpha^{\prime}(G)=\lfloor\frac{n}{4}\rfloor+1$ and $s=t$, by Lemma \ref{Lemma 2.6} (2), $uv\in E(\overline{G})$ for each vertex $u\in V_{12} $ and $v \in V_{2}$. Since $|V_{22}|\geq 3$, we may assign number $1$ to $s-2$ independent edges in $[V_{12}\setminus \{v_{1}, v_{2}\}, V_{21}]$ in $\overline{G}$, and $\frac{1}{2}$ to each edge of a cycle with length $n-2t-(s-2)$ in $\overline{G}[V_{2}]$. For a fractional matching $\bar{f}$ of $\overline{G}$, let $\bar{f}(v_{1}v_{2})=1$ in $\overline{G}$. Thus we have
\begin{align*}
\alpha^{\prime}(\overline{G}) \geq 1+s-2+\frac{n-2t-(s-2)}{2} =\frac{n-s}{2}=\frac{n-t}{2}.
\end{align*} \par
When $s<t$ and there exists a $\frac{1}{2}-$cycle in $V_{12}\setminus X$, it can be proved similar to the above in (1). Assume that there is no $\frac{1}{2}-$cycle in $V_{12}\setminus X$ and there are $p$ edges assigned number $1$ in $E[V_{12}]$, where $p\geq 1$ and $p$ is an integer. Since $f(v)=1$ for every $v\in V_{1}$, $p=0$ is equivalent to $s=t$ or $G[V_{12}\setminus X]$ is a collection of disjoint $\frac{1}{2}-$cycles, which has been discussed above. If $p=1$, that is $t=s+1$ and $|V_{12}|=2t-s=s+2$, then there exactly exists an edge say $ww_{1}$ assigned number $1$ in $G[V_{12}\setminus X ]$. Without loss of generality, there exists a vertex $x$ in $V_{21}$ such that $xw\in E(\overline{G})$ by Lemma \ref{Lemma 2.5} (1). Since $|V_{22}|\geq 3$, let $\bar{f}(v_{1}v_{2})=\bar{f}(xw)=\bar{f}(y_{1}w_{1})=1$ in $\overline{G}$, $y_{1} \in V_{22}$. Then for each vertex $v\in V_{12}\setminus \{v_{1}, v_{2}, w, w_{1}\}$ and each vertex $u\in V_{2}\setminus\{x, y_{1}\}$, we have $uv\in E(\overline{G})$. We may assign number $1$ to $s-2$ independent edges in $[V_{12}\setminus \{v_{1}, v_{2}, w, w_{1}\}, V_{21}\setminus \{x\}]$, and $\frac{1}{2}$ to each edge of a cycle with length $n-2t-s$ in $\overline{G}[V_{2}]$. Thus we have
\begin{align*}
\alpha^{\prime}(\overline{G})\geq 1+1+1+(s-2)+\frac{n-2t-s}{2}=\frac{n-s}{2}>\frac{n-t}{2}.
\end{align*} \par
If $p\geq 2$, then there exist two edges $w_{1}w_{2}$ and $w_{3}w_{4}$ assigned number $1$ in $G[V_{12}\setminus X]$. Without loss of generality, there exist two vertices $x_{1}, x_{2}$ in $V_{21}$ such that $x_{1}w_{1}, x_{2}w_{3} \in E(\overline{G} )$ and let $\bar{f}(x_{1}w_{1})=\bar{f}(x_{2}w_{3})=1$ in $\overline{G} $. 
In $\overline{G}$, there are $s-2$ independent edges in $[X\setminus \{v_{1}, v_{2}\}, V_{21}\setminus \{x_{1}, x_{2}\}]$ by Lemma \ref{Lemma 2.6} (2) and we may assign number $1$ to them. Noting that $n-2t-s\geq 2t-2s-2$, we may assign number $1$ to $2t-2s-2$ independent edges in $\overline{G}[V_{12}\setminus\{X\cup \{w_{1}, w_{3}\}\}, V_{22}]$. If 
$n-2t-(2t-s-2)=0$, that is $n=4t-s-2$ and $|V_{22}|=2t-2s-2$, we may assign number $1$ to edge $v_{1}v_{2}$ in $\overline{G}$. Thus we have $\alpha^{\prime}(\overline{G}) \geq 1+2t-s-2=2t-s-1=\frac{n-s}{2}>\frac{n-t}{2}$.\par
If $n-2t-(2t-s-2)=1$, that is $n=4t-s-1$ and $|V_{22}|=2t-2s-1$, then we may choose $y \in V_{22}$ and assign number $\frac{1}{2}$ to edges of $C_{3}$ induced by $v_{1}, v_{2}$ and $y$ in $\overline{G}$, and $1$ to $2t-2s-2$ independent edges in $\overline{G}[V_{12}\setminus\{X\cup \{w_{1}, w_{3}\}\}, V_{22}\setminus\{y\}]$. Thus we have $\alpha^{\prime}(\overline{G}) \geq \frac{3}{2}+2t-s-2=2t-s-\frac{1}{2}=\frac{n-s}{2}>\frac{n-t}{2}$.\par

If $n-2t-(2t-s-2)=2$, that is $n=4t-s$ and $|V_{22}|=2t-2s$, then we may choose $y_{1}, y_{2} \in V_{22}$ and assign number $1$ to edges $v_{1}v_{2}$ and $y_{1}y_{2}$ in $\overline{G}$, and $1$ to $2t-2s-2$ independent edges in $\overline{G}[V_{12}\setminus\{X\cup \{w_{1}, w_{3}\}\}, V_{22}\setminus\{y_{1}, y_{2}\}]$. Thus we have $\alpha^{\prime}(\overline{G}) \geq 1+2t-s-2+1=2t-s=\frac{n-s}{2}>\frac{n-t}{2}$.\par
When $n-2t-(2t-s-2)\geq 3$, we may assign number $\frac{1}{2}$ to each edge of a cycle with length $n-2t-(2t-s-2)$ in $\overline{G}[V_{2}]$, and $1$ to edge $v_{1}v_{2}$ in $\overline{G}$. Thus we have
\begin{align*}
\alpha^{\prime}(\overline{G}) \geq 1+2t-s-2+\frac{n-2t-(2t-s-2)}{2} =\frac{n-s}{2}> \frac{n-t}{2}.
\end{align*}\par
Therefore, $\alpha^{\prime}(\overline{G}) \geq \frac{n-t}{2}\geq \frac{n}{4}+3$ for $n\geq 28$. This completes the proof.       $\hfill\square $\\
\vspace{8pt}

Without loss of generality, we may assume that $\alpha^{\prime}(G) \leq \alpha^{\prime}(\overline{G}) $ as follows.\\
\begin{theorem}\label{Theorem 3.3}
	Let $G$ be a graph of order $n\geq 28$. If both $G$ and $\overline{G}$ are not empty, then 
	\begin{align*}
	\alpha^{\prime}(G)+\alpha^{\prime}(\overline{G})\geq\frac{n+1}{2}
	\end{align*}
	with equality holds if and only if $G\cong K_{1, n-1}$.
\end{theorem}
\noindent\textbf{Proof} Since $G$ is not empty, $\alpha^{\prime}(G)\geq1$. By Lemma \ref{Lemma 2.11} (1)-(5), when $1\leq \alpha^{\prime}(G)<3$, $\alpha^{\prime}(G)+\alpha^{\prime}(\overline{G})\geq \frac{n+1}{2}$. The equality holds if and only if $G\cong K_{1, n-1}$ by Lemma \ref{Lemma 2.11} (1) and Lemma \ref{Lemma 2.8} (1).\par 
When $3\leq \alpha^{\prime}(G)\leq \frac{n}{4}$, by Lemma \ref{Lemma 3.1} (1), we have 
$$\alpha^{\prime}(G)+\alpha^{\prime}(\overline{G})\geq\alpha^{\prime}(G)+\frac{n-\alpha^{\prime}(G)}{2}=\frac{n+\alpha^{\prime}(G)}{2}\geq \frac{n+3}{2}>\frac{n+1}{2}.$$Then we only need to consider the case that $\frac{n}{4}<\alpha^{\prime}(G)\leq \frac{n}{2}$. \par
If $n\equiv0, 1(mod4)$, then $\frac{n+1}{4}\leq\lfloor\frac{n}{4} \rfloor+\frac{1}{2}\leq \alpha^{\prime}(G)\leq \frac{n}{2}$. So we have 
\begin{align*}
\alpha^{\prime}(G)+\alpha^{\prime}(\overline{G})\geq 2\alpha^{\prime}(G)\geq 2(\lfloor \frac{n}{4}\rfloor +\frac{1}{2})\geq \frac{n+1}{2}.
\end{align*}
The equality holds if and only if $n\equiv1(mod4)$ and $\alpha^{\prime}(G)=\alpha^{\prime}(\overline{G})=\lfloor\frac{n}{4} \rfloor +\frac{1}{2}$, which is a contradiction to Lemma \ref{Lemma 3.2} (1) in virtue of $\lfloor\frac{n}{4} \rfloor +\frac{1}{2}<\frac{n}{4}+3$.               \par      
If $n\equiv2, 3(mod4)$, $\frac{n}{4}<\lfloor\frac{n}{4} \rfloor+1\leq \alpha^{\prime}(G)\leq \frac{n}{2}$.
When $\alpha^{\prime}(G)=\lfloor\frac{n}{4}\rfloor +1$, by Lemma \ref{Lemma 3.2} (2), $\alpha^{\prime}(\overline{G})\geq \frac{n}{4}+3$. Thus $\alpha^{\prime}(G)+\alpha^{\prime}(\overline{G})\geq \lfloor\frac{n}{4}\rfloor +1+\frac{n}{4}+3\geq \frac{n}{2}+3>\frac{n+1}{2}$. 
When $\frac{n+2}{4}<\lfloor\frac{n}{4}\rfloor +\frac{3}{2}\leq \alpha^{\prime}(G) \leq \frac{n}{2}$, we have
\begin{align*}
\alpha^{\prime}(G)+\alpha^{\prime}(\overline{G})\geq 2\alpha^{\prime}(G)\geq 2(\lfloor \frac{n}{4}\rfloor +\frac{3}{2})>\frac{n+2}{2}>\frac{n+1}{2}.
\end{align*}\par

Therefore, $\alpha^{\prime}(G)+\alpha^{\prime}(\overline{G})>\frac{n+1}{2}$. This completes the proof. $\hfill\square $\par
The following result is implied by Theorem \ref{Theorem 3.3}.
\begin{lemma}\label{Theorem 3.4}
      For a graph of order $n$, $n\geq 2$, we have
      \begin{align*}
      \alpha^{\prime}(G)+\alpha^{\prime}(\overline{G})\geq \frac{n}{2}
      \end{align*}
      with equality holds if and only if $G$ is an empty graph or $G$ is a complete graph.
\end{lemma}
The following result is deduced by Lemma \ref{Lemma 2.9}, Lemma \ref{Lemma 2.10} and Lemma \ref{Lemma 2.11}.
\begin{lemma}\label{Lemma 2.12}
	If $\alpha^{\prime}(G)=2$ or $\alpha^{\prime}(G)=\frac{5}{2}$, $\overline{G}$ contains no isolated vertices, then $\alpha^{\prime}(\overline{G}) =\frac{n}{2}$ for $n\geq 10$.
\end{lemma}

\begin{theorem}\label{Theorem 3.5}
	 If both $G$ and $\overline{G}$ contain no isolated vertices and $n\geq 28$. Then 
	\begin{align*}
	\alpha^{\prime}(G)+\alpha^{\prime}(\overline{G})\geq\frac{n+4}{2}
	\end{align*}
	and the equality holds if and only if  $G \cong K_{2}(p, q; \ell)$ or $K_{1, m} \cup K_{1, n-2-m}\subseteq G \subseteq K_{2, n-2}$, where $p, q, \ell, m$ are non-negative integers, $q\geq 1$ and $1\leq p, m\leq n-3$.          
\end{theorem}
\noindent\textbf{Proof}  Since both $G$ and $\overline{G}$ contain no isolated vertices, by Lemma \ref{Lemma 2.8}, we have $\alpha^{\prime}(G)$ and $ \alpha^{\prime}(\overline{G})\geq 2$. When $s=0$, since $G$ contains no isolated vertices, $V_{2}$ is empty and $G[V_{1}]$ contains a fractional perfect matching by Lemma \ref{Lemma 2.4}. The assumption $\alpha^{\prime}(G)\leq \alpha^{\prime}(\overline{G})$ implies $\alpha^{\prime}(\overline{G})= \alpha^{\prime}(G)=\frac{n}{2}$. Thus $\alpha^{\prime}(G)+\alpha^{\prime}(\overline{G})=\frac{n}{2} +\frac{n}{2}=n>\frac{n+4}{2}$.\par
When $\alpha^{\prime}(G)=2$ or $\frac{5}{2}$, by Lemma \ref{Lemma 2.12}, $\alpha^{\prime}(\overline{G})=\frac{n}{2}$. It follows that $\alpha^{\prime}(G)+\alpha^{\prime}(\overline{G})\geq\frac{n+4}{2}$. The equality holds when  $\alpha^{\prime}(G)=2$.\par 
When $\alpha^{\prime}(G)=3$,
if $s=1$ or $2$, we have $\alpha^{\prime}(\overline{G})\geq \frac{n-s+1}{2}$ by Lemma \ref{Lemma 3.1} (2). It follows that $\alpha^{\prime}(G)+\alpha^{\prime}(\overline{G}) \geq 3+\frac{n-s+1}{2}>\frac{n+4}{2}$. 
If $s=t=3$, we have $\alpha^{\prime}(\overline{G})\geq \frac{n-3}{2}+1=\frac{n-1}{2}$ by Lemma \ref{Lemma 3.1} (3), and then $\alpha^{\prime}(G)+\alpha^{\prime}(\overline{G}) \geq 3+\frac{n-1}{2}=\frac{n+5}{2}>\frac{n+4}{2}$. \par
When $\alpha^{\prime}(G)=\frac{7}{2}$, then $1\leq s\leq 3$. By Lemma \ref{Lemma 3.1} (2), we have $\alpha^{\prime}(\overline{G})\geq \frac{n-s+1}{2}$. It follows that $\alpha^{\prime}(G)+\alpha^{\prime}(\overline{G}) \geq \frac{7}{2}+\frac{n-s+1}{2}\geq \frac{7}{2}+\frac{n-3+1}{2}=\frac{n+5}{2}>\frac{n+4}{2}$.\par
When $4\leq\alpha^{\prime}(G)=t\leq \frac{n}{4}$, since $s\geq 1$, $\alpha^{\prime}(\overline{G})\geq \frac{n-s+1}{2}$ by Lemma \ref{Lemma 3.1} (2). Thus
\begin{align*}
\alpha^{\prime}(G)+\alpha^{\prime}(\overline{G}) \geq t+\frac{n-s+1}{2}\geq t+\frac{n-t+1}{2}=\frac{n+t+1}{2}\geq \frac{n+5}{2}>\frac{n+4}{2}.
\end{align*}
Then we only need to consider the case that $\frac{n}{4}<\alpha^{\prime}(G)\leq \frac{n}{2}\leq \alpha^{\prime}(\overline{G})$.\par
If $n\equiv0, 1(mod4)$, then $\lfloor\frac{n}{4} \rfloor+\frac{1}{2}\leq \alpha^{\prime}(G)\leq \frac{n}{2}$. When $\alpha^{\prime}(G)=\lfloor \frac{n}{4} \rfloor+\frac{1}{2}$ or $\lfloor \frac{n}{4} \rfloor+1$, by Lemma \ref{Lemma 3.2}, $\alpha^{\prime}(\overline{G})\geq\frac{n}{4}+3$. Thus
\begin{align*}
\alpha^{\prime}(G)+\alpha^{\prime}(\overline{G}) \geq \lfloor \frac{n}{4} \rfloor+\frac{1}{2}+\frac{n}{4}+3>\frac{n}{2}+3>\frac{n+4}{2}.
\end{align*}
When $\alpha^{\prime}(G)\geq \lfloor\frac{n}{4}\rfloor +\frac{3}{2}>\frac{n+4}{4}$, we obtain $\alpha^{\prime}(G)+\alpha^{\prime}(\overline{G})\geq 2\alpha^{\prime}(G) \geq 2\lfloor\frac{n}{4}\rfloor +3 >\frac{n+4}{2}$.\par
If $n\equiv2, 3(mod4)$, $\lfloor\frac{n}{4} \rfloor+1\leq \alpha^{\prime}(G)\leq \frac{n}{2}$. When $\alpha^{\prime}(G)=\lfloor\frac{n}{4}\rfloor +1$ or $\lfloor\frac{n}{4}\rfloor +\frac{3}{2}$, by Lemma \ref{Lemma 3.2}, $\alpha^{\prime}(\overline{G})\geq \frac{n}{4}+3$. Thus we have
$\alpha^{\prime}(G)+\alpha^{\prime}(\overline{G})\geq\lfloor\frac{n}{4}\rfloor +1+\frac{n}{4}+3> \frac{n}{2}+3>\frac{n+4}{2}$. When $\alpha^{\prime}(G)\geq \lfloor\frac{n}{4}\rfloor +2> \frac{n+4}{4}$, we have
\begin{align*}
\alpha^{\prime}(G)+\alpha^{\prime}(\overline{G})\geq 2\alpha^{\prime}(G)\geq 2(\lfloor \frac{n}{4}\rfloor +2)>\frac{n+4}{2}.
\end{align*}\par

Therefore, $\alpha^{\prime}(G)+\alpha^{\prime}(\overline{G})>\frac{n+4}{2}$.\par
If $\alpha^{\prime}(G)+\alpha^{\prime}(\overline{G})=\frac{n+4}{2}$, then $\alpha^{\prime}(G)=2$. By Lemma \ref{Lemma 2.9} and both $G$ and $\overline{G}$ contain no isolated vertices, we have $G\cong K_{2}(p,q;\ell)$ or $K_{1, m} \cup K_{1, n-2-m}\subseteq G \subseteq K_{2, n-2}$, where $p, q, \ell, m$ are non-negative integers, $q\geq 1$ and $1\leq p, m\leq n-3$. Therefore, we complete the proof. $\hfill\square $\\



\end{document}